\documentclass[12pt]{article}
\usepackage{amsmath}
\usepackage{amssymb}
\usepackage{vatola}

\def\q{\quad}

\def\t{\hbox}
\def\mod#1{\ (\hbox{\rm mod}\ #1)}
\def\qtq#1{\q\t{#1}\q}
\def\f{\frac}
\def\e{\equiv}
\def\b{\binom}

\def\Ls#1#2{\Big(\f{#1}{#2}\Big)}
\let \pro=\proclaim
\let \endpro=\endproclaim
\begin{document}

 \centerline {\bf
Inequalities for binomial coefficients} \par\q\newline
\centerline{Zhi-Hong Sun} \par\q\newline \centerline{School of
Mathematical Sciences, Huaiyin Normal University,}
\centerline{Huaian, Jiangsu 223001, P.R. China} \centerline{E-mail:
zhihongsun@yahoo.com} \centerline{Homepage:
http://www.hytc.edu.cn/xsjl/szh}

\abstract{In this paper we prove several inequalities for binomial
coefficients. For instance, if $ k$ and $n$ are positive integers
such that $n\ge 400$ and $[\frac n5]\le k\le [\frac n2]$, where
$[x]$ is the greatest integer not exceeding $x$, then
$$\binom nk<\Big(1-\frac{5(k-[\f n5])}{6n^2}\Big) \frac{n^{n-\f
12}}{k^k(n-k)^{n-k}}.$$
\par\q
\newline MSC: Primary 05A20; Secondary 11B65, 33B15.
 \newline Keywords:
 binomial coefficient; factorial; inequality.}
 \endabstract
 \footnotetext[1] {The author is
supported by the Natural Sciences Foundation of China (grant no.
11371163).}

\section*{1. Introduction}
\par\q Since the binomial coefficients are quotients of factorials of
positive integers, in order to determine the size of binomial
coefficients it is very natural to use good approximations for
factorials of positive integers. The famous Stirling's approximation
[1] asserts that
$$n!\sim \sqrt{2\pi n}\Big(\f ne\Big)^n\qtq{as} n\rightarrow
+\infty.$$ More precisely,  Robbins [6] and Feller [2] proved that
for any positive integer $n$,
$$
\sqrt{2\pi n } \left(  \frac{n}{e}\right)^n e^{\frac{1}{12n+1}} < n!
< \sqrt{2 \pi n }  \left(  \frac{n}{e}\right)^n e^{\frac{1}{12n}}.
 \tag 1.1$$
 Hence, for $m,n\in\Bbb N$ with $m\ge 2$,
 $$\align\b{mn}n=\f{(mn)!}{n!((m-1)n)!}&<\f{{\sqrt{2\pi mn}}(\f
 {mn}{\t{e}})^{mn}\t{e}^{\f 1{12mn}}}
 {{\sqrt{2\pi n}}(\f
 {n}{\t{e}})^{n}\t{e}^{\f 1{12n+1}}{\sqrt{2\pi (m-1)n}}(\f
 {(m-1)n}{\t{e}})^{(m-1)n}\t{e}^{\f 1{12(m-1)n+1}}}
\\&=\sqrt{\f m{2\pi
(m-1)n}}\Ls{m^m}{(m-1)^{m-1}}^n\t{e}^{\f 1{12mn}-\f 1{12n+1}-\f
1{12(m-1)n+1}}.\endalign $$ As
$$\align&\f 1{12mn}-\f 1{12n+1}-\f
1{12(m-1)n+1}\\&=\f{-144mn-12n+12m+1-144mn^2(m-1)-12mn}{12mn(12n+1)
(12(m-1)n+1)}<0,\endalign$$ we obtain
$$\b {mn}n<\sqrt{\f m{2\pi
(m-1)n}}\Ls{m^m}{(m-1)^{m-1}}^n\qtq{for $m,n\in\Bbb N$ and $m\ge
2$.}\tag 1.2$$
 \par For any real number $x>0$ let
$$ \Gamma(x) = \int_{0}^{ + \infty} t^{x-1} e^{-t} dt.$$
It is well known that $\Gamma(n+1)=n!$ for any positive integer $n$.
 A conjecture of Srinivasa Ramanujan [5] asserts that for
$x>0$, $$\Gamma(x+1 ) =\sqrt{\pi} \left( \frac{x}{e}\right)^x \left(
8x^3 +4x^2 +x +\frac{\theta(x)}{30} \right)^\frac{1}{6} ,$$
 where $\theta (x) \to 1$ as $x\to + \infty$ and
  $\frac{3}{10} < \theta (x) <1$. When $x=n\ge 2$ is a positive integer,
  in 2006 Hirschhorn [3] established
the following better result:
$$ n!=\sqrt {\pi} \left( \frac{n}{e}\right)^n \left( 8n^3+4n^2+n+\frac{1}{30}-\frac{11}{240n}+
 \frac{r}{240n^2} \right)^{\frac{ 1}{6}},\tag 1.3$$
 where $5<r<11$. Note that $\frac{
r}{240n}-\frac{11}{240}<\frac{11}{240n}-\frac{11}{240}<0$. From
(1.3) we have
$$\sqrt {\pi}\left(
\frac{n}{e}\right)^n (8n^3+4n^2+n)^{\frac{ 1}{6}} <n!<\sqrt
{\pi}\left( \frac{n}{e}\right)^n \left( 8n^3+4n^2+n+\frac{1}{30}
\right)^{\frac{ 1}{6}}\ \t{for}\ n\ge 2 .\tag 1.4$$
 Let $\Bbb N$ be
the set of positive integers. In 2001 St$\check {\rm a}$nic$\check
{\rm a}$ [7] showed that for any $m,n,p\in\Bbb N$ with $m>p$,
$$\t{e}^{-\f 1{8n}}\sqrt{\f m{2\pi p(m-p)n}}\Big(\f
{m^m}{(m-p)^{m-p}}\Big)^n<\b{mn}{pn}<\sqrt{\f m{2\pi
p(m-p)n}}\Big(\f {m^m}{(m-p)^{m-p}}\Big)^n.\tag 1.5$$
\par Let $[x]$
 be the greatest integer not exceeding $x$. Suppose $k,n\in\Bbb N$
 with $2\le k\le [\f n2]$. In this paper we improve
 the known inequality $\b
nk\le\f{n^n}{k^k(n-k)^{n-k}}$ $(1\le k<n)$ ([4]) by proving that
$$\b nk<\Big(1-\f{5(k-1)}{6n^2}\Big)
\f{n(n-1)^{n-1}}{k^k(n-k)^{n-k}}\tag 1.6$$ and
$$\b nk<\Big(1-\f{5(k-[\f n5])}{6n^2}\Big) \f{n^{n-\f
12}}{k^k(n-k)^{n-k}}\ \t{for}\ k\ge \big[\f n5\big]\ \t{and}\ n\ge
400.\tag 1.7$$
\par
 Suppose $m,n\in\Bbb N$ with $m,n\ge 2$.
  In this paper we  use (1.2) to establish
the following inequality:
$$ \b n{[\f nm]}<\f m{\sqrt{2\pi
(m-1)n}}\Ls{m}{(m-1)^{\f{m-1}m}}^n\q \t{for}\q  n\ge 2m-1.\tag 1.8$$
If $k,n\in\Bbb N$, $n\ge 6$, $k\le \f n2$ and $\f{2^n}{n+1}<\b nk\le
\f{2^n}n$, we also show that $k$ is the least positive integer $r$
such that $\b nr>\f{2^n}{n+1}$.

\section*{2. Main results}
\pro{Lemma 2.1} Let $n\in\Bbb N$, $n\ge 3$ and $k\in\{0,1,\ldots,
[\f n2]-1\}$. Then
$$\b nkk^k(n-k)^{n-k}>\t{\rm e}^{\f {11}{12n^2}}\b n{k+1}(k+1)^{k+1}(n-k-1)^{n-k-1}.$$
\endpro
Proof. Set
$$F(n,k)=\f{\b nkk^k(n-k)^{n-k}}{\b
n{k+1}(k+1)^{k+1}(n-k-1)^{n-k-1}}.$$ For $k=0$ we see that
$$F(n,k)=\f{n^n}{n(n-1)^{n-1}}=(1+\f 1{n-1})^{n-1}>
2>\sqrt{\t{e}}>\t{e}^{\f {11}{12n^2}}.$$ Thus the result is true for
$k=0$. Now we assume that $1\le k\le \f n2-1$. It is clear that
$$\b n{k+1}=\f{n(n-1)\cdots(n-k+1)(n-k)}{(k+1)!}=\f{n-k}{k+1}\b
nk.$$ Thus,
$$ F(n,k)
=\f{k^k(n-k)^{n-k-1}}{(k+1)^k(n-k-1)^{n-k-1}}=\Big(1-\f
1{k+1}\Big)^k\Big(1+\f 1{n-k-1}\Big)^{n-k-1}$$ and so
$$\log F(n,k)=k\log\Big(1-\f 1{k+1}\Big)+(n-k-1)\log\Big(1+\f
1{n-k-1}\Big).$$ For $m\ge 1$ we see that
$$\f 1{nm^n}>\f 1{(n+1)m^{n+1}}\qtq{and}
\lim_{n\rightarrow+\infty}\f 1{nm^n}=0.$$ Thus,
$$\log\big(1+\f 1m\big)=\f 1m-\f 1{2m^2}+\f 1{3m^3}-\f
1{4m^4}+\cdots >\f 1m-\f 1{2m^2}+\f 1{3m^3}-\f 1{4m^4}$$
 and so
$$m\log\big(1+\f 1m\big)>1-\f 1{2m}+\f 1{3m^2}
-\f 1{4m^3}\q\t{for $m\ge 1$}.\tag 2.1$$ As $1\le k\le \f n2-1$ we
have $\f n2\le  n-k-1\le n-2$. Hence, taking $m=n-k-1$ in (2.1) we
get
$$\aligned (n-k-1)\log\big(1+\f 1{n-k-1}\big)&>1-\f 1{2(n-k-1)}
+\f 1{3(n-k-1)^2}-\f 1{4(n-k-1)^3}\\&>1-\f 1n+\f 1{3(n-2)^2}-\f
2{n^3}.\endaligned\tag 2.2$$ On the other hand, for $m\ge 2$ we have
$$\align -\log\big(1-\f 1m\big)&=\f 1m+\f 1{2m^2}+\f 1{3m^3}+\f
1{4m^4}+\f 1{5m^5}+\f 1{6m^6}+\cdots
\\&< \f 1m+\f 1{2m^2}+\f 1{3m^3}+\f
1{4m^4}+\f 1{4m^5}+\f 1{4m^6}+\cdots
\\&=\f 1m+\f 1{2m^2}+\f 1{3m^3}+\f
1{4(m-1)m^3}=\f{12m^3-6m^2-2m-1}{12(m-1)m^3}\endalign$$ and so
$$(m-1)\log\big(1-\f 1m\big)>-1+\f 1{2m}
+\f 1{6m^2}+\f 1{12m^3}\q\t{for $m\ge 2$}.\tag
2.3$$ For $1\le k\le \f n2-1$ we have $k+1\le \f n2$. Taking $m=k+1$
in (2.3) we get
$$k\log\big(1-\f 1{k+1}\big)>-1+\f 1{2(k+1)}+\f 1{6(k+1)^2}
+\f 1{12(k+1)^3}>-1+\f 1n+\f 2{3n^2}+\f 2{3n^3}.\tag 2.4$$ Using
 (2.2) and (2.4) we obtain
$$\align \log F(n,k)
&=k\log\big(1-\f 1{k+1}\big)+(n-k-1)\log\big(1+\f 1{n-k-1}\big)
\\&>-1+\f 1n+\f 2{3n^2}+\f 2{3n^3}+1-\f 1n+\f 1{3(n-2)^2}-\f 2{n^3}
\\&=\f {11}{12n^2}+\f 13\Big(\f 1{(n-2)^2}-\f 3{4n^2}-\f 4{n^3}\Big)
\\&=\f {11}{12n^2}+\f{n^3-4n^2+52n-64}{12n^3(n-2)^2}
>\f {11}{12n^2},\endalign$$
Thus, $F(n,k)> \t{e}^{\f {11}{12n^2}}.$ This proves the lemma.
\pro{Theorem 2.1} Let $n\in\Bbb N$, $n\ge 4$ and
$k\in\{2,3,\ldots,[\f n2]\}$. Then
$$\b nk<\t{e}^{-\f {11(k-1)}{12n^2}}\f{n(n-1)^{n-1}}{k^k(n-k)^{n-k}}
<\Big(1-\f{5(k-1)}{6n^2}\Big) \f{n(n-1)^{n-1}}{k^k(n-k)^{n-k}}.$$
\endpro
Proof. For $m\in\{0,1,\ldots,[\f n2]\}$ set $g(m)=\b
nmm^m(n-m)^{n-m}$. By Lemma 2.1,
$$\f{g(m)}{g(m-1)}<\t{e}^{-\f {11}{12n^2}}\qtq{for}m=1,2,\ldots,\big[\f n2\big].$$
Thus,
$$g(k)=\f{g(k)}{g(k-1)}\cdot\f{g(k-1)}{g(k-2)}
\cdots\f{g(2)}{g(1)}\cdot g(1)<\t{e}^{-(k-1)\f
{11}{12n^2}}n(n-1)^{n-1}.$$ For $0<x<1$ we have
$\t{e}^{-x}=1-x+\f{x^2}{2!}-\f{x^3}{3!}+\cdots<1-x+\f{x^2}2$. Thus,

$$\aligned \t{e}^{-(k-1)\f
{11}{12n^2}}&<1-\f{11(k-1)}{12n^2}+\f{121(k-1)^2}{288n^4}
\\&=1-\f{5(k-1)}{6n^2}-\f{k-1}{12n^2}\Big(1-\f{121(k-1)}{24n^2}\Big)
<1-\f {5(k-1)}{6n^2}.\endaligned\tag 2.5$$
Now combining all the
above we prove the theorem.

 \pro{Theorem 2.2} For any positive integer $n$ and
 $k\in\{0,1,2,\ldots,n\}$ we have
 $$\b nk<\sqrt{\f 2{\pi}}\cdot \f{2^n}{\sqrt n}<\f
 45\cdot\f{2^n}{\sqrt n}.$$
 \endpro
Proof. Clearly the result is true for $n=1$. Now assume that $n\ge
2$. When $n$ is even, taking $m=2$ and substituting $n$ with $n/2$
in (1.2) we see that
$$\b n{\f n2}<\sqrt{\f 2{\pi}}\cdot \f{2^n}{\sqrt n}.$$
When $n$ is odd, from the above we deduce that
$$\b n{[\f n2]}=\b n{\f{n-1}2}=\f{2n}{n+1}\b{n-1}{\f{n-1}2}
<\f{2n}{n+1}\sqrt{\f 2{\pi}}\cdot \f{2^{n-1}}{\sqrt {n-1}}.$$ As
$n^3<n^3+n^2-n-1=(n+1)^2(n-1)$ we have $n\sqrt n<(n+1)\sqrt{n-1}$.
Thus
$$\f n{(n+1)\sqrt{n-1}}<\f 1{\sqrt n}\qtq{and so}
\b n{[\f n2]}<\sqrt{\f 2{\pi}}\cdot\f n{(n+1)\sqrt{n-1}}\cdot 2^n
<\sqrt{\f 2{\pi}}\cdot \f{2^n}{\sqrt n}$$ as asserted. To complete
the proof, we note that $\b nk\le \b n{[n/2]}$ and $\sqrt {\f
2{\pi}}<\f 45$.

\pro{Theorem 2.3} Let $m,n\in\Bbb N$ with $m,n\ge 3$ and $n\ge
2m-1.$ Then
$$\b n{[\f nm]}<\f m{\sqrt{2\pi
(m-1)n}}\Ls{m}{(m-1)^{\f{m-1}m}}^n.$$
\endpro
Proof.  Assume $n\e r\mod m$ with $r\in\{0,1,\ldots,m-1\}$. Then
$[\f nm]=\f{n-r}m$.  For $r=0$ the result follows from (1.2). Now we
assume $1\le r\le m-1$. It is clear that
$$\b nr\b{n-r}{\f{n-r}m}=\b n{n-r}\b{n-r}{\f{n-r}m}=\b
n{\f{n-r}m}\b{n-\f{n-r}m}{n-r-\f{n-r}m}=\b
n{\f{n-r}m}\b{n-\f{n-r}m}r$$ and so
$$\b n{\f{n-r}m}=\f{\b nr}{\b{n-\f{n-r}m}r}\b{n-r}{\f{n-r}m}
=\b{n-r}{\f{n-r}m}\prod_{k=0}^{r-1}\f{n-k}{n-\f{n-r}m-k}.$$ For
$k\in\{0,1,\ldots,r-1\}$ we have
$$\f{(m-1)(n-k+1)}m\ge n-\f{n-r}m-k=\f{(m-1)n-(km-r)}m
>\f{(m-1)(n-k)}m$$ and so
$$\f m{m-1}\cdot\f{n-k}{n-k+1}\le \f{n-k}{n-\f{n-r}m-k}<\f m{m-1}.$$
Hence,
$$\Ls m{m-1}^r\f{n+1-r}{n+1}\le
\prod_{k=0}^{r-1}\f{n-k}{n-\f{n-r}m-k}<\Ls m{m-1}^r$$
 and therefore
$$\b{n-r}{\f{n-r}m}\Ls m{m-1}^r\f{n+1-r}{n+1}\le \b n{\f{n-r}m}
<\b{n-r}{\f{n-r}m}\Ls m{m-1}^r.\tag 2.6$$
 From (2.6) and (1.2) we deduce that
$$\align \b n{\f{n-r}m}&<\Ls m{m-1}^r\b{n-r}{\f{n-r}m}
<\Ls m{m-1}^r\cdot \f m{\sqrt{2\pi
(m-1)(n-r)}}\Ls{m^m}{(m-1)^{m-1}}^{\f {n-r}m}
\\&=\sqrt{\f n{n-r}}(m-1)^{-\f rm}\cdot
\f m{\sqrt{2\pi (m-1)n}}\Ls{m^m}{(m-1)^{m-1}}^{\f nm}.\endalign$$
Set $x=(m-1)^{\f 2{m}}$. Then $x>1$. For $r\ge 2$ we see that
$$r< x^{r-1}+x^{r-2}+\cdots+x+1=\f{x^r-1}{x-1}\qtq{and so}
\f r{x^r-1}< \f 1{x-1}.$$ Since $(m-1)^2\ge 4>\t{e}>(1+\f 1m)^m$ we
see that $(m-1)^{\f 2m}>1+\f 1m$ and so $m>\f 1{(m-1)^{\f 2m}-1}$.
As $1\le r\le m-1$, we see that
$$n\ge 2m-1 > m-1+\f 1{(m-1)^{\f 2{m}}-1}>r+\f r{(m-1)^{\f {2r}{m}}-1}.$$
Thus,
$$\f{n-r}r>\f 1{(m-1)^{\f {2r}{m}}-1}\qtq{and so}\f n{n-r}=1+\f
r{n-r}<(m-1)^{\f{2r}{m}}.$$ Therefore
$$\sqrt{\f n{n-r}}<(m-1)^{\f rm}\qtq{and so}
 \b n{\f{n-r}m}<\f m{\sqrt{2\pi
(m-1)n}}\Ls{m^m}{(m-1)^{m-1}}^{\f nm}.$$ This proves the theorem.

\pro{Lemma 2.2} Let $k\in\Bbb N$, $k\ge 80$ and $r\in\{0,1,2,3,4\}$.
Then
$$\b {5k+r}k<\f{(5k+r)^{5k+r-\f 12}}
{k^k(4k+r)^{4k+r}}.$$
\endpro
Proof. Using Maple we know that the result is true for $80\le k\le
200$. Now assume $k>200$. For $r=0$, from (1.2) we have
$$\b{5k}k<\sqrt{\f 5{8\pi k}}\cdot\Ls{5^5}{4^4}^k=\sqrt{\f{25}{8\pi}}\cdot
\f {(5k)^{5k}}{\sqrt{5k}\cdot k^k(4k)^{4k}}<\f
{(5k)^{5k}}{\sqrt{5k}\cdot k^k(4k)^{4k}} .$$ So the result is true
for $r=0$.  Now suppose $r=1$. Clearly
$$\align \b{5k+1}k=\f{5k+1}{4k+1}\b{5k}k&<\f{5k+1}{4k+1}\cdot
\sqrt{\f{25}{8\pi}}\cdot \f {(5k)^{5k}}{\sqrt{5k}\cdot
k^k(4k)^{4k}}\\&=\sqrt{\f{25}{8\pi}}\cdot\f{(\f {5k}{5k+1})^{5k-\f
12}}{(\f {4k}{4k+1})^{4k}}\cdot \f{(5k+1)^{5k+1-\f
12}}{k^k(4k+1)^{4k+1}}.\endalign$$ By Lagrange's mean value theorem,
for $m\ge 1$ there exists a real number $\theta\in (0,1)$ such that
$$\log (1+\f 1m)=\f{\log(m+1)-\log m}{m+1-m}=\f 1{m+\theta}.$$  Thus, $\f 1{m+1}<\log(1+\f 1m)<\f 1m$ for $m\ge 1$.
Using this inequality we see that
$$\align &\Big(5k-\f 12\Big)\log\f{5k}{5k+1}-4k\log\f{4k}{4k+1}
\\&=-\Big(5k-\f 12\Big)\log\Big(1+\f 1{5k}\Big)+4k\log(1+\f 1{4k})
\\&<-\Big(5k-\f 12\Big)\f 1{5k+1}+4k\cdot \f 1{4k}=\f 3{2(5k+1)}.
\endalign$$
Hence,
$$\f{(\f {5k}{5k+1})^{5k-\f
12}}{(\f {4k}{4k+1})^{4k}}<\t{e}^{\f 3{2(5k+1)}}\q \t{and so}
\q\b{5k+1}k<\sqrt{\f{25}{8\pi}}\cdot\t{e}^{\f 3{2(5k+1)}}\cdot
\f{(5k+1)^{5k+1-\f 12}}{k^k(4k+1)^{4k+1}}.$$ For $k\ge 114$ we have
$$\sqrt{\f{25}{8\pi}}\cdot\t{e}^{\f 3{2(5k+1)}}<1\qtq{and so}
\b{5k+1}k<\f{(5k+1)^{5k+1-\f 12}}{k^k(4k+1)^{4k+1}}.$$ This proves
the result in the case $r=1$.
\par Now assume $r\in\{2,3,4\}$.
By (1.4),
$$\align &(5k+r)!<\sqrt
{\pi}\Ls{5k+r}{\t{e}}^{5k+r}\Big(8(5k+r)^3+4(5k+r)^2+5k+r+\f
1{30}\Big)^{\f 16},
\\&k!>\sqrt{\pi}\Ls k{\t{e}}^k(8k^3+4k^2+k)^{\f 16},
\\&(4k+r)!>\sqrt{\pi}\Ls{4k+r}{\t{e}}^{4k+r}(8(4k+r)^3
+4(4k+r)^2+4k+r)^{\f 16}.\endalign$$ Thus,
$$\align
&\b{5k+r}k=\f{(5k+r)!}{k!(4k+r)!}
\\&<\f{(5k+r)^{5k+r}}{\sqrt{\pi}k^k(4k+r)^{4k+r}}
\Ls{8(5k+r)^3+4(5k+r)^2+5k+r+\f 1{30}}{(8k^3+4k^2+k)(8(4k+r)^3
+4(4k+r)^2+4k+r)}^{\f 16}.\endalign$$ To prove the result, it is
sufficient to show that
$$\f{8(5k+r)^3+4(5k+r)^2+5k+r+\f 1{30}}{(8k^3+4k^2+k)(8(4k+r)^3
+4(4k+r)^2+4k+r)}<\f{\pi^3}{(5k+r)^3}.\tag 2.7$$ Set
$$\aligned F(k,r)&=\pi^3(8k^3+4k^2+k)(8(4k+r)^3
+4(4k+r)^2+4k+r)\\&\q-(5k+r)^3\Big(8(5k+r)^3+4(5k+r)^2+5k+r+\f
1{30}\Big).\endaligned$$ Then
$$\align F(k,r)&>3.14^3(8k^3+4k^2+k)(8(4k+2)^3
+4(4k+2)^2+4k+2)\\&\q-(5k+4)^3(8(5k+4)^3+4(5k+4)^2+5k+4+1)
\\&>1808k^6-343032k^5-1019794k^4-1260084k^3-810552k^2-270342k-37185.\endalign$$
As $k>200$ we have
$$\align &37185<k^2<k^6,\ 270342k<k^4<k^6,\ 810552k^2<k^5<k^6,
\\&1260084k^3<k^6,\ 1019794k^4<26k^6,\ 343032k^5<1750k^6.\endalign$$
Thus,
$$F(k,r)>(1806-1750-26-1-1-1-1)k^6=26k^2>0$$
and so (2.7) holds. Hence the result is true for $r=2,3,4$. Now the
proof is complete.
 \pro{Lemma 2.3} Let $n\in\Bbb N$ and
$k\in\{0,1,\ldots,n\}$. Then
$$ \f{n^n}{k^k(n-k)^{n-k}}\le 2^n.$$
\endpro
Proof. Clearly the result is true for $k=0,n$. Now assume $1\le
k<n$.  For $1\le x\le \frac{n}{2}$, let
$$f(x)=x^x (n-x)^{n-x}\qtq{and}
 F(x)=\log f(x)=x\log
x+(n-x)\log(n-x).$$ Then
$$F'(x)=\log x-\log(n-x)\quad
\text{and} \quad F^{\prime\prime}(x)=\frac{1}{x} +
\frac{1}{n-x}=\frac{n}{x(n-x)}>0.$$
  By the symmetry, we
may assume $1\le k\le \f n2$. Then $F^{\prime} (k) < F^{\prime}
\left(\frac{n}{2} \right)=0$ and
  so $f(k)\geq
f(\frac{n}{2})$. Thus $$\frac{n^n}{k^k (n-k)^{n-k}}= \frac{n^n}{
f(k)} \leq  \frac{n^n}{ f(\frac{n}{2})} =2^n.$$ The proof is now
complete.

 \pro{Theorem 2.4} Let
$k,n\in\Bbb N$, $n\ge 400$ and $[\f n5]\le k\le [\f n2]$. Then
$$\align \b nk&<\t{e}^{-(k-[\f n5])\f {11}{12n^2}}
\f{n^{n-\f 12}}{k^k(n-k)^{n-k}}<\Big(1-\f{5(k-[\f n5])}{6n^2}\Big)
\f{n^{n-\f 12}}{k^k(n-k)^{n-k}} \\&<\Big(1-\f{5(k-[\f
n5])}{6n^2}\Big)\f{2^n}{\sqrt n}.\endalign$$
\endpro
Proof. By Lemmas 2.1, 2.2 and (2.5),
$$\align \b nkk^k(n-k)^{n-k}&=\b n{[\f n5]}\big[\f n5\big]^{[\f n5]}
\Big(n-\big[\f n5\big]\Big)^{n-[\f n5]} \prod_{s=[\f n5]}^{k-1}\f{\b
n{s+1}(s+1)^{s+1}(n-s-1)^{n-s-1}}{\b nss^s(n-s)^{n-s}}
\\&<n^{n-\f 12}\t{e}^{-(k-[\f n5])\f {11}{12n^2}}<n^{n-\f 12}
\Big(1-\f{5(k-[\f n5])}{6n^2}\Big).
\endalign$$ This together with Lemma 2.3 gives the result.

\pro{Corollary 2.1} Let $k,n\in\Bbb N$, $n\ge 400$ and $[\f n5]\le
k\le [\f {4n}5]$. Then
$$\b nk<\f{n^{n-\f 12}}{k^k(n-k)^{n-k}}.$$
\endpro
Proof. By the symmetry, we may assume $[\f n5]\le k\le [\f n2]$. Now
the result follows from Theorem 2.4 immediately.
\par{\bf Remark 2.1} Corollary 2.1 was conjectured by
the author's student Shu-Yao Yi.
\par Let $n\in\{3,4,5,\ldots\}$. Following [4] we define $f(n)$ to be the
least positive integer $k$ such that $\b nk>\f{2^n}{n+1}$.
 \pro{Theorem 2.5} Suppose
$k,n\in\Bbb N$, $n\ge 6$, $k\le \f n2$ and $\f{2^n}{n+1}<\b nk\le
\f{2^n}n$. Then $k=f(n)$.\endpro
 Proof. As $\binom nr \leq \binom{n}{
[\frac{n}{2}] }$ for $r=0,1,\ldots,n$, we see that
$$(n-3)\binom{n}{[\frac{n}{2}]}\ge \sum_{r=2}^{n-2}\binom nr
 =2^n-2-2n.$$
Thus,
$$\align\b n{[\f n2]-1}&=\f{[n/2]}{n+1-[n/2]}\b n{[\f n2]}\ge
\f{2^n-2-2n}{n-3}\cdot\f{2[n/2]}{2n+2-2[n/2]}
\\&\ge \f{2^n-2-2n}{n-3}\cdot\f{n-1}{n+3}=\f{2^n}{n+1}\cdot
\f{2^n-2-2n}{2^n}\cdot\f{n^2-1}{n^2-9}.\endalign$$ Clearly
$2^{7+1}>7^2(7-2)$. If $2^{m+1}>m^2(m-2)$ for $m\ge 7$, then
$$2^{m+2}>2m^2(m-2)=(m+1)^2(m+1-2)+m^3-5m^2+m+1>(m+1)^2(m+1-2).$$
Thus, by induction we have $2^{n+2}>(n+1)^2(n+1-2)$ for $n\ge 6$.
Hence $\f{n+1}{2^{n-1}}<\f 8{n^2-1}$ and so
$$\f{2^n-2-2n}{2^n}=1-\f{n+1}{2^{n-1}}>1-\f
8{n^2-1}=\f{n^2-9}{n^2-1}.$$ Now, from the above we deduce that
$$\b
n{[\f n2]-1}>\f{2^n}{n+1}\qtq{and so} f(n)\le [\f n2]-1. \tag 2.8$$
Therefore, $$f(n)\le \f{n-2}2<\f{n^2-n-1}{2n+1}\qtq{and so}
\f{n-f(n)}{f(n)+1}>\f{n+1}n.$$
 By the definition,
$\b n{f(n)}>\f{2^n}{n+1}$. Thus, from the above we deduce that
$$\b n{f(n)+1}=\f{n-f(n)}{f(n)+1}\b n{f(n)}>\f{n-f(n)}{f(n)+1}\cdot
\f{2^n}{n+1}>\f{2^n}n.$$ Hence $k=f(n)$. This proves the theorem.
\par\q
\par By doing calculations with the help of Maple, we pose the
following conjecture.

 \pro{Conjecture 2.1} There are infinitely many pairs $(n,k)$
 $(k,n\in\Bbb N, k\le \f n2)$  such that
 $\f{2^n}{n+1}<\b nk\le \f{2^n}n$.\endpro
 \par The first few examples of $(n,k)$ $(n\le 1500)$ are listed below:
$$\align &(2,1),(4,1),(19,6),(61,23),(89,35),(130,53),(139,57),
(291,126),(343,150),(521,233),
\\&(712,323),(788,359),(929,426),(950,436),(971,446),(1080,498),
(1289,598),(1387,645).\endalign$$
\par By Theorem 2.5, Conjecture 2.1 is equivalent to the following conjecture.
\pro{Conjecture 2.2} There are infinitely many positive integers $n$
such that $\b n{f(n)}\le \f{2^n}n$.\endpro
 \pro{Conjecture 2.3} We
have
$$\lim_{x\rightarrow+\infty}\f{|\{n\le x:\ \b n{f(n)}<\f{2^n}n\}|}
{\sqrt { x/\log x}}=\root 3\of{\f{\pi}2}.$$
\endpro


\begin{thebibliography}{99}
\bibitem{Bateman} H. Bateman, Higher Transcendental Functions,
  vol. I, McGraw-Hill, New York, 1953, p.47.

\bibitem{Fell} W. Feller, An Introduction to Probability Theory and Its
Applications, 3${}^{\text{rd}}$ edition, New York, Wiley, 1968,
pp.52-54.
\bibitem{14}M.D. Hirschhorn, A new version of Stirling's formula, Math.
 Gazette 90(2006), 286-292.

 \bibitem{KKS} D. Kim, A. Sankaranarayanan and Z.H. Sun,
On the properties of a sequence originated from
   binomial coefficients, arXiv:1309.7567.

\bibitem{Ram} S. Ramanujan, The lost notebook and other unpublished
   papers (ed. by S. Raghavan and
  S. S. Rangachari), Springer, 1988.

  \bibitem{Robin} H. Robbins, Remarks of Stirling's formula,
   Amer. Math. Monthly {\bf 62}(1955), 26-29.

 \bibitem{Stanica} P. St$\check {\rm a}$nic$\check {\rm a}$,
 Good lower and upper bounds on
 binomial coefficients, J. Inequal. Pure and Appl. Math. {\bf 3}
 (2001), Art.30, 5pp.


\end{thebibliography}
\end{document}